\documentclass[11pt]{article}
\usepackage{graphicx}
\usepackage{tikz}
\usetikzlibrary{arrows.meta}
\usepackage{setspace}
\usepackage{booktabs}
\usepackage{capt-of}
\usepackage{etoolbox}
\usepackage{hyperref}

\BeforeBeginEnvironment{tabular}{\begin{singlespace}}
\AfterEndEnvironment{tabular}{\end{singlespace}}

\usepackage[margin=0.6in]{geometry}

\begin{document}

\title {A Trivalent Graph of Girth 14 and Order 372}

\author{Geoffrey Exoo\\
\small Department of Mathematics and Computer Science\\[-0.8ex]
\small Indiana State University\\[-0.8ex]
\small Terre Haute, IN 47809, USA\\[-0.8ex]
\small\tt geoffrey.exoo@gmail.com\\
\\
Theodore Kolokolnikov\\
\small Department of Mathematics and Statistics\\[-0.8ex]
\small Dalhousie University\\[-0.8ex]
\small Halifax, Nova Scotia, Canada\\[-0.8ex]
\small\tt tkolokol@gmail.com\\
\\
}

\maketitle

\begin{abstract}
The smallest known trivalent graph of girth 14
is constructed.  The graph has 372 vertices, improving
the old upper bound of 384 \cite{exoo14}.
The construction uses methods that may be new to this problem.
\end{abstract}

\begin{center}
AMS Subject Classifications: {05C25, 05C35}\hfill
\end{center}

\vspace{5mm}

\section{Introduction}

This paper deals with the cage problem, which asks for the identification of the smallest graphs
with given degree and girth \cite{exjaj}.  A {\em $(k,g)$-graph} is a $k$-regular graph with girth $g$.
A {\em $(k,g)$-cage} is a $(k,g)$-graph with a minimum possible order.
In this paper we construct the smallest known $(3,14)$-graph.  Our graph has order 372.

The construction uses ideas that appeared in our paper on an 
extremal problem on algebraic connectivity \cite{ekjs}.
There, we considered extremal regular graphs of given diameter that had the
minimum possible second eigenvalue (i.e. largest algebraic connectivity),
and showed that such a graph of diameter 13 has to have order 380
and must consist of
two cubic Moore trees each having 190 vertices, with additional edges 
connecting the leafs between two trees (see Figure~\ref{fig380}, top).
One of the graphs constructed in this paper, in addition to being
a girth-14 graph of order 380, also has the maximum possible
algebraic connectivity ($0.351227$) among cubic graphs with diameter 13.
By excising 8 vertices from this new graph, we obtain
a new record $(3,14)$ graph of order 372.

\section{Construction of the (3,14)-Graph of order 372}

We begin by reviewing some notation and terminology that is used
throughout the paper. We then describe the construction of a
(3,14) graph of order 380
which we excise to obtain a (3,14) record of order 372.

\subsection{Moore Trees}

A {\em Moore tree} of degree $k$ and radius $r$ is a tree with a central vertex such that every
leaf is at distance $r$ from the central vertex and every non-leaf vertex has degree $k$.
In this paper we consider {\em cubic Moore trees}, i.e., the case $k=3$.

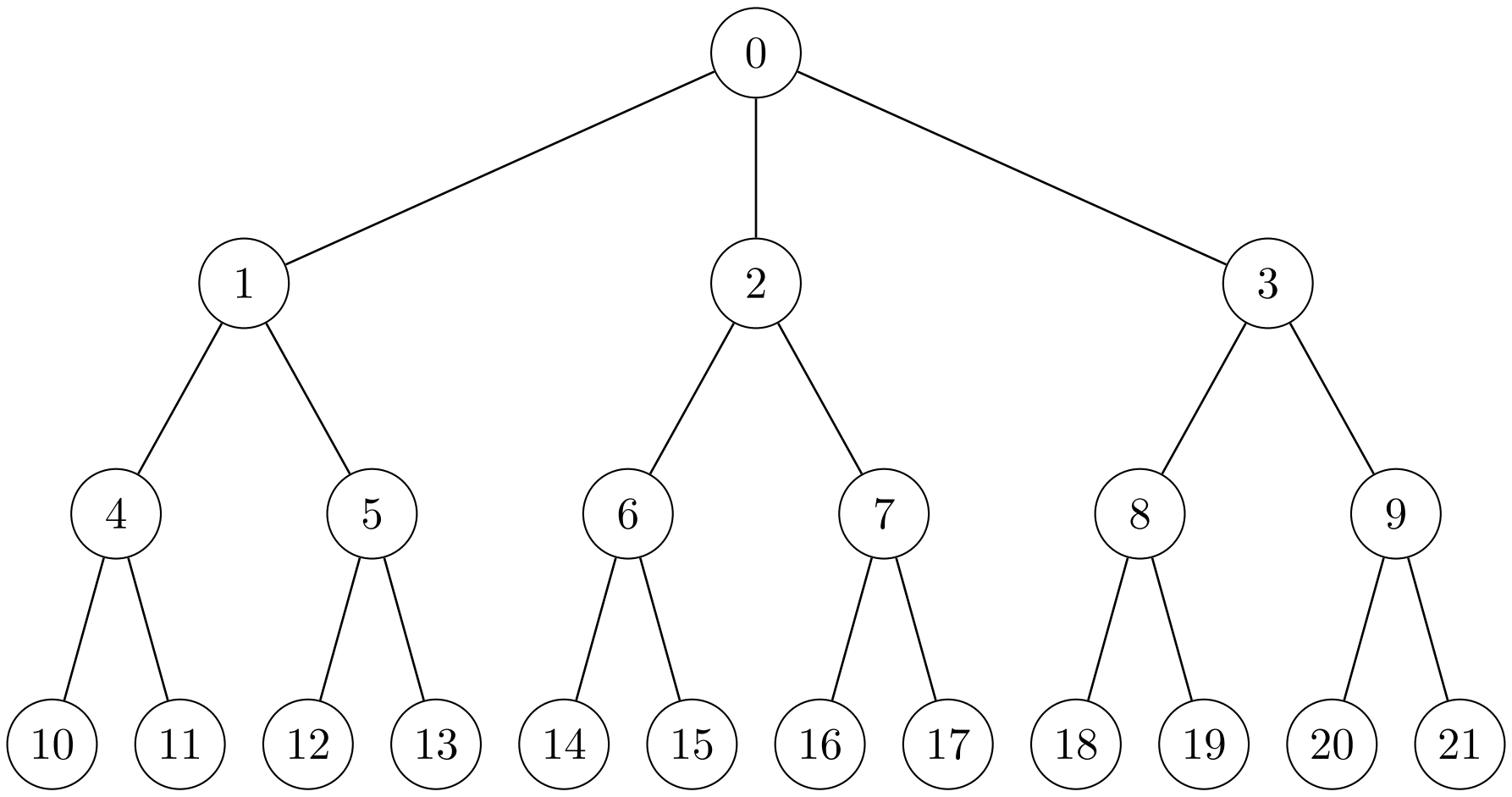
\begin{figure}
\center
\begin{tikzpicture}[
    vertex/.style={circle, draw, minimum size=7mm, inner sep=0pt, font=\small},
    edge/.style={line width=0.5pt}
  ]
  \def\radius{3}
  \def\xsep{1}
  \def\ysep{1.8}

  \pgfmathtruncatemacro{\L}{3*2^(\radius-1)}
  \node[vertex] (v0) at (0,0) {0};
  \foreach \d in {1,...,\radius} {
    \pgfmathtruncatemacro{\n}{3*2^(\d-1)}
    \pgfmathtruncatemacro{\last}{\n-1}
    \pgfmathtruncatemacro{\w}{2^(\radius-\d)}
    \pgfmathtruncatemacro{\off}{\n-2}
    \foreach \k in {0,...,\last} {
      \pgfmathsetmacro{\x}{\xsep*(\w*(\k+0.5) - \L/2)}
      \pgfmathtruncatemacro{\v}{\off+\k}
      \node[vertex] (v\v) at (\x,-\d*\ysep) {\v};
    }
  }

  \foreach \v in {1,2,3} {
    \draw[edge] (v0) -- (v\v);
  }
  \foreach \d in {2,...,\radius} {
    \pgfmathtruncatemacro{\n}{3*2^(\d-1)}
    \pgfmathtruncatemacro{\last}{\n-1}
    \pgfmathtruncatemacro{\off}{\n-2}
    \pgfmathtruncatemacro{\poff}{\n/2-2}
    \foreach \k in {0,...,\last} {
      \pgfmathtruncatemacro{\v}{\off+\k}
      \pgfmathtruncatemacro{\p}{\poff+floor(\k/2)}
      \draw[edge] (v\p) -- (v\v);
    }
  }
\end{tikzpicture}
\caption{A cubic Moore tree of radius 3 with the standard labeling.}
\label{treefig}
\end{figure}

The {\em standard labeling} of the vertices of a cubic Moore Tree is given as follows.
The vertices of the tree are labeled  $u_0, \dots, u_{3 ( 2^r - 1)}$.
For $r > 0$,
vertex $u_0$ is the central vertex and is adjacent to vertices $u_1$, $u_2$, and $u_3$.
The remaining edges in the tree are of the form
$(u_i,u_{2i+2})$ and $(u_i,u_{2i+3})$,
for $1 \leq i \leq 3 (2^{r-1} - 1)$.
We shall refer to $u_1$, $u_2$, and $u_3$ as the
children of $u_0$ and to $u_{2i+2}$ and $u_{2i+3}$ as the
children of $u_i$.
The standard labeling for a cubic Moore tree of radius 3
is shown in Figure~\ref{treefig}.

The {\em level sets} of the cubic Moore tree are the sets of vertices at a given distance from the
central vertex.  For the tree shown in Figure 1, there are four level sets:
$L_0$, $L_1$, $L_2$, and $L_3$.  Level set 3 is the set of leaves of the tree in the figure, or in general,
level set $r$ is the set of leaves.

\subsection{Constructing a (3,14)-Graph of order 380}

\begin{figure}[!tbp]
\center
\input{fig_sigma2.tex}
\caption{A graphical depiction of tree automorphism $\sigma$.}
\label{tree6fig}

\vspace{\floatsep}
{\small
\setlength{\tabcolsep}{3pt}
\begin{minipage}[t]{0.71\textwidth}
\centering
\begin{tabular}{|c|c|c|l|}
\hline
Cycle Length & Level & Labels & Cycle \\ \hline
  1 & 0 && $(0)$ \\
  3 & 1 && $(1,3,2)$ \\
  3 & 2 && $(4,8,6)$ \\
  3 & 2 && $(5,9,7)$ \\
  6 & 3 && $(10,19,15,11,18,14)$ \\
  6 & 3 && $(12,21,17,13,20,16)$ \\
 12 & 4 && $(22,41,33,25,39,31,23,40,32,24,38,30)$ \\
 12 & 4 && $(26,45,37,29,43,35,27,44,36,28,42,34)$ \\
 12 & 5 && $(46,84,68,52,80,64,48,82,66,50,78,62)$ \\
 12 & 5 && $(47,85,69,53,81,65,49,83,67,51,79,63)$ \\
 12 & 5 && $(54,92,76,60,88,72,56,90,74,58,86,70)$ \\
 12 & 5 && $(55,93,77,61,89,73,57,91,75,59,87,71)$ \\
 12 & 6 &$A_0,B_0$& $(94,170,138,106,162,130,98,166,134,102,158,126)$ \\
 12 & 6 &$A_1,B_1$& $(95,171,139,107,163,131,99,167,135,103,159,127)$ \\
 12 & 6 &$A_2,B_2$& $(96,172,140,108,164,132,100,168,136,104,160,128)$ \\
 12 & 6 &$A_3,B_3$& $(97,173,141,109,165,133,101,169,137,105,161,129)$ \\
 12 & 6 &$A_4,B_4$& $(110,186,154,122,178,146,114,182,150,118,174,142)$ \\
 12 & 6 &$A_5,B_5$& $(111,187,155,123,179,147,115,183,151,119,175,143)$ \\
 12 & 6 &$A_6,B_6$& $(112,188,156,124,180,148,116,184,152,120,176,144)$ \\
 12 & 6 &$A_7,B_7$& $(113,189,157,125,181,149,117,185,153,121,177,145)$ \\ \hline
\end{tabular}
\captionof{table}{Disjoint Cycles in the Cycle Decomposition of $\sigma$ and $\tau$.
These permutations are fully determined by the vertex rotation at the root and the flips at $u_4$, $u_5$, $u_{10}$ and $u_{12}$, as shown in Figure~\ref{tree6fig}.}
\label{cycletab}
\end{minipage}\hfill
\begin{minipage}[t]{0.26\textwidth}
\centering
\begin{tabular}{|c|c|r|}
\hline
A-Orbit & B-Orbit & power of $\sigma$ \\ \hline
$A_0$ & $B_0$ &  6 \\
$A_0$ & $B_0$ &  1 \\
$A_1$ & $B_4$ &  1 \\
$A_1$ & $B_7$ &  0 \\
$A_2$ & $B_3$ &  3 \\
$A_2$ & $B_3$ &  4 \\
$A_3$ & $B_5$ &  3 \\
$A_3$ & $B_6$ & 10 \\
$A_4$ & $B_1$ &  9 \\
$A_4$ & $B_7$ &  7 \\
$A_5$ & $B_2$ &  7 \\
$A_5$ & $B_6$ &  9 \\
$A_6$ & $B_2$ &  0 \\
$A_6$ & $B_5$ &  4 \\
$A_7$ & $B_1$ & 10 \\
$A_7$ & $B_4$ &  6 \\ \hline
\end{tabular}
\captionof{table}{Edges joining the two trees.}
\label{edgetab}
\end{minipage}%
}
\end{figure}

Our first construction is of a $(3,14)$-graph of order 380.
We were able to complete this construction finding three different
graphs.  We present the construction of the graph of order 380 that led us
to the $(3,14)$-graph of order 372.

We start with two copies of a cubic Moore tree of radius 6, $T_1$ and $T_2$.
Each such tree has order 190 and is given the standard labeling.
The vertices of $T_1$ are labeled  $u_0, \dots, u_{189}$, and those of $T_2$ are labeled $v_0, \dots, v_{189}$.

The graph is completed by adding edges from
each of the leaf vertices in $T_1$ (vertices $u_{94}, \dots, u_{189}$) to
two leaf vertices in $T_2$.

The number of ways to do this is too large for a simple exhaustive search,
so the search space is reduced by forcing symmetry.
To do this automorphisms $\sigma$ of $T_1$ and $\tau$ of $T_2$ are chosen.
At first we simplified things
by making $\sigma$ and $\tau$ the same automorphism, i.e., defining $\tau$ such that
$\tau(v_i) = v_j$ whenever $\sigma(u_i) = u_j$.  Two of the three $(3,14)$-graphs of order 380
that were found used the same automorphism
for $\sigma$ and $\tau$. The graph of order 380 that led us to the graph of order 372 was
one of these two.

So we continue the construction assuming that $\sigma$ and $\tau$ are the same automorphism,
acting on the two different trees.

For heuristic considerations, we were looking for tree automorphisms that had only
one fixed point (the central vertex).
Such an automorphism $\sigma$ of a cubic Moore tree with the standard labeling can be defined as follows.

\begin{itemize}
\item[1.] $\sigma(u_0) = u_0$.
\item[2.] $\sigma$ cyclically permutes $u_1$, $u_2$, and $u_3$.  Either $3$-cycle can be used.
This guarantees that $u_0$ is the only fixed point of $\sigma$.
\item[3.] For each $i$, $ 1 \leq i \leq 3(2^{r-1}-1)$, $\sigma$ either fixes the children of $u_i$
or interchanges them.
\end{itemize}

We speculated (on admittedly very little evidence) that $\sigma$ should be chosen so that
the order of $\sigma$ (and the size of the leaf orbits) should be something close
to the square root of the number of leaves.  After much experimentation (and testing of automorphisms
of larger orders) we found that order 12 might be the best choice.

Upon further experimentation, we finally found the $\sigma$ that works. 
It is constructed by cyclically permuting the children of vertex $u_0$, and then
swapping the two subtrees of vertices $u_{4}$, $u_{5}$, $u_{10}$ and $u_{12}$.
Figure \ref{tree6fig} illustrates this construction. 
Our choice of $\sigma$ is given in Table~\ref{cycletab}, where each row in the list gives
one disjoint cycle.  $\sigma$ is the product of all $20$ of these cycles.
Next note that $|L_6|$, the number of leaves in one of our trees, is $96$.
Therefore $\sigma$ partitions $L_6$ into $|L_6|/12 = 8$ classes.

The last eight cycles give the partition of the leaves of $T_1$ (and of $T_2$) into
orbits of size $12$.
To specify the edges joining the leaves of the two trees we use an idea that is
reminiscent of a voltage construction.

First we label the eight leaf orbits of $\sigma$ by $A_0, \ldots, A_7$ as shown in Table~\ref{cycletab}.
We similarly label leaf orbits of $\tau$ by $B_0, \ldots, B_7$.
The entries for $A$ are the subscripts of $u_i$ in $T_1$ and similarly for $B$;
for example $A_0=(u_{94}, u_{170}, \ldots, u_{126})$ and $B_0=(v_{94}, v_{170}, \ldots, v_{126})$.


In Table~\ref{edgetab} the connections between the leaves of
$T_1$ and the leaves of $T_2$ are given.  Each row determines a matching
between one of the A-sets (column 1) and one of the B-sets (column 2).
The third column gives the power of $\sigma$
that is applied to the B-sets.  Denote the first
(leftmost or minimum) element in each of the $A_i$ by $r_i$.  So
$r_0 = 94, r_1 = 95, \dots , r_7 = 113$.  Let  {\tt $A_i$   $B_j$   $k$}
be a row from the table.
This row determines $12$ edges $(u,v)$ that are added to the
graph.
The edges are $(u,v) = (\sigma^t(u_{r_i}), \tau^{t+k}(v_{r_j}))$, for $0 \leq t < 12$.
Table~\ref{edgetab} together with $\sigma, \tau$ (obtained via four flips and one rotation) 
fully determine the resulting graph. 

The final step, verifying that the resulting cubic graph has girth $14$,
was done using the {\tt pickg} program in B. McKay's {\em nauty} package \cite{nauty}.
The entire graph is shown in Figure \ref{fig380}(top).
This completes the construction of a $(3,14)$-graph of order $380$.

\subsection{Using Excision on the graph of order 380 to obtain the (3,14)-Graph of order 372}

The last step in the construction of the $(3,14)$-graph of order $372$
is an excision.  First we remove all eight of the level $0$ and $1$ vertices from
the cubic Moore trees $T_1$ and $T_2$, i.e., vertices
$u_0$, $u_1$, $u_2$ and $u_3$ from $T_1$, and vertices
$v_{0}$, $v_{1}$, $v_{2}$ and $v_{3}$ from $T_2$.
We then add edges in such a way as to obtain a girth-14 cubic graph.
Of the 10395 ways to pair up the 12 vertices of degree 2, exactly one gives girth 14.

%

\begin{figure}[!tbp]
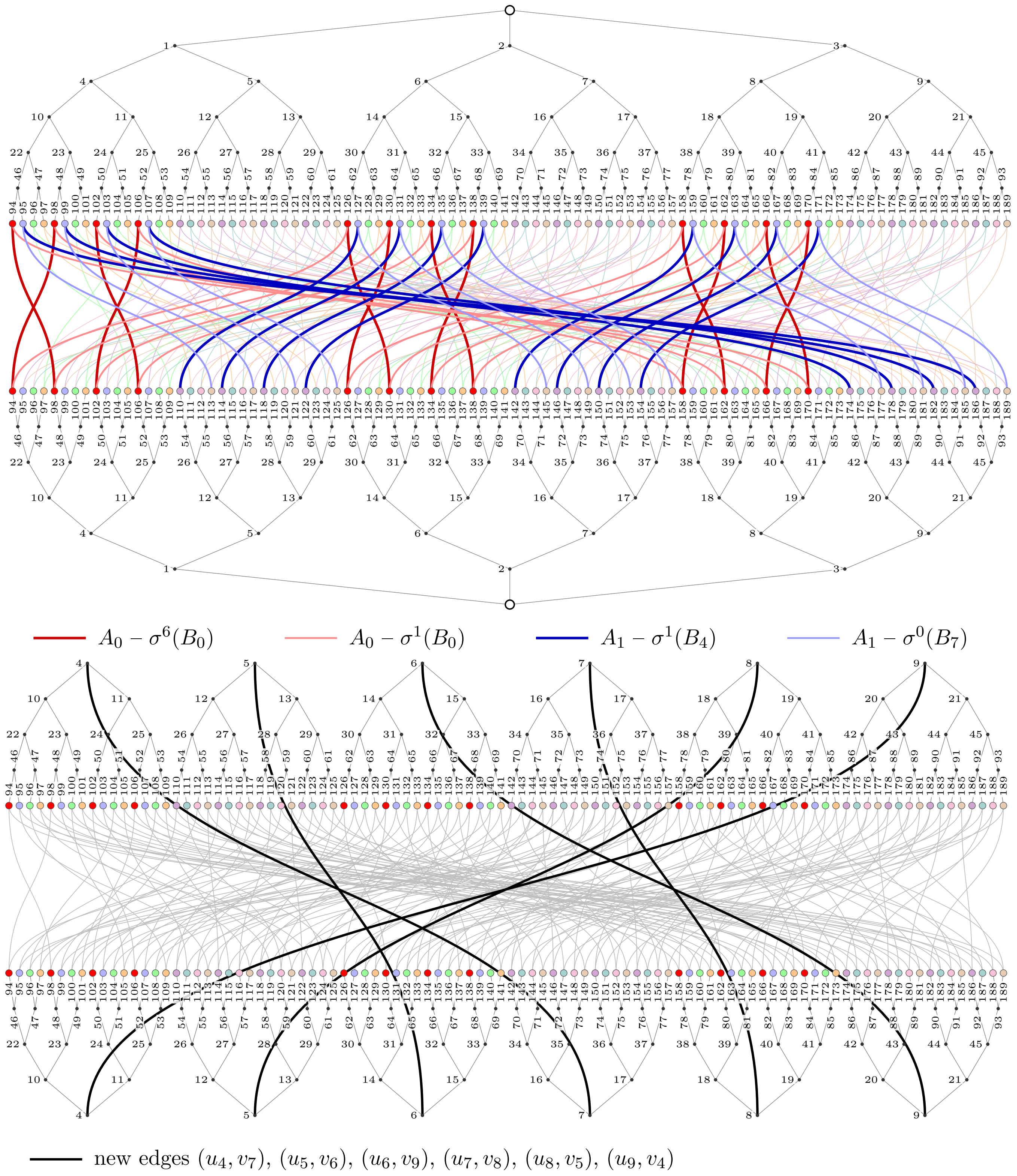

\center
\input{fig_full380.tex}
\input{fig_excise372.tex}
\caption{Top: the $(3,14)$-graph of order 380. 
Leafs are colored by which $A$ or $B$-set they belong as in Figure~\ref{tree6fig}. 
Edges between the trees are joined according to Table~\ref{edgetab}, 
coloured according to the $A$-set they originate from.
Bottom: the $(3,14)$-graph of order 372 obtained via excision from the 380 graph.
The roots and its children were deleted from both trees, then the remaining vertices 
were joined in a unique way that forms a girth-14 cubic graph.
}
\label{fig380}
\end{figure}

The resulting graph is $3$-regular and has order $372$, girth 14, and is shown in Figure \ref{fig380}(bottom); the added edges are listed in the caption.
It can be verified, again using {\em nauty}, that the girth of the graph is $14$. 
In the appendix we list three programs, the first shows how
tree automorphisms are generated.
The second
program uses a candidate permutation for $\sigma$ and outputs the edge list for the
graph constructed using $\sigma$ as above.
The third program does the excision that yields the graph of order $372$.

The record 372-vertex graph of girth 14 can be downloaded from the
\href{https://houseofgraphs.org/meta-directory/cages}{The House of Graphs} ~\cite{hog}.
It has diameter 12 and 18 antipodal vertex pairs. 
Its automorphism group is $(S_4\times D_8):C_2$ (SmallGroup(384,18036) in GAP) of size 384. It consists of six edge Moore trees joined to each-other. 
Its vertex count is 1.465 times that of the Moore bound of 254 (the smallest theoretically possible number for (3,14)-cage). 
Its second eigenvalue is $\lambda_2 = 2.640252$ so that
its algebraic connectivity is 3-2.640252=0.359748, below the upper bound of 0.371404 for diameter-12 graph \cite{ekjs}. This is in contrast
to the 380-vertex graph it was excised from, which has diameter 13 and attains the maximum possible algebraic connectivity of 0.351227 for its diameter.

After completion of this paper, we were able to reconstruct the graph using a very different method.
This method will be part of an upcoming paper \cite{ekhunt} focusing on 
extremal questions relating
diameter and algebraic connectivity.

\section{Acknowledgements}

The $(3,14)$-graph of order $372$ which is the main topic of this paper
can be generated using the three programs shown in the Appendix.
These programs, along with equivalent (but much faster) C programs,
were written by the authors of this paper with no AI assistance.
The LLM Claude (Anthropic) was used in two places in our work.
\begin{itemize}
\item[1.] A search procedure for the automorphism $\sigma$ was generated by Claude.  An equivalent
search was designed without Claude's assistance.  Both searches were successful, in that they
both found the graph of order $380$ that was used in the final construction.  
\item[2.] The figures were created by Claude.
\end{itemize}

\newpage

\appendix

\section{Program 1}

The first program generates a cubic Moore tree automorphism.  This program
uses the list of swaps that generates our graph of order 380.  In the search phase,
the program generates random lists of from 3 to 7 randomly selected nodes.
The random selection was biased toward low number nodes in a variety of ways.

\begin{verbatim}

# A Python program that converts a recipe for a cubic Moore Tree automorphism into
# a permutation.  The swaplist variable is the recipe.

N = 190       # order of cubic Moore tree
S =  94       # number of non-leaf nodes
L =  12       # size of group
T =   4       # number of swaps (random choice in range 3 to 7)

swaplist = [4, 5, 10, 12]    # in general a random list of T nodes from 1 to 93

def swap(sig, u, v, n):

    if u >= N:
        return
    for i in range(n):
        sig[u+i], sig[v+i] = sig[v+i], sig[u+i]
    swap(sig, 2*u+2, 2*v+2, 2*n)

def buildperm():

    sig = list(range(N))
    swap(sig, 2, 3, 1)    # the 3-cycle as a composition of two transpositions
    swap(sig, 1, 2, 1)

    for i in range(T):
        node = swaplist[i]
        lchild = 2*node + 2
        rchild = 2*node + 3
        swap(sig, lchild, rchild, 1)
    return sig

if __name__ == "__main__":

    sigma = buildperm()
    print(sigma)

\end{verbatim}

\section{Program 2}

The second program uses the tree automorphism $\sigma$ and prints the edgelist of the graph,
along with the orbits of the group generated by $\sigma$.

\begin{verbatim}

#
# A Python program that generates the graph of order 380 described in this paper.
#

N = 380
H = 190        # Half of N

# An automorphism of the cubic Moore tree of radius 6.

g0 = [
   0,    3,    1,    2,    8,    9,    4,    5,    6,    7,
  19,   18,   21,   20,   10,   11,   12,   13,   14,   15,
  16,   17,   41,   40,   38,   39,   45,   44,   42,   43,
  22,   23,   24,   25,   26,   27,   28,   29,   30,   31,
  32,   33,   34,   35,   36,   37,   84,   85,   82,   83,
  78,   79,   80,   81,   92,   93,   90,   91,   86,   87,
  88,   89,   46,   47,   48,   49,   50,   51,   52,   53,
  54,   55,   56,   57,   58,   59,   60,   61,   62,   63,
  64,   65,   66,   67,   68,   69,   70,   71,   72,   73,
  74,   75,   76,   77,  170,  171,  172,  173,  166,  167,
 168,  169,  158,  159,  160,  161,  162,  163,  164,  165,
 186,  187,  188,  189,  182,  183,  184,  185,  174,  175,
 176,  177,  178,  179,  180,  181,   94,   95,   96,   97,
  98,   99,  100,  101,  102,  103,  104,  105,  106,  107,
 108,  109,  110,  111,  112,  113,  114,  115,  116,  117,
 118,  119,  120,  121,  122,  123,  124,  125,  126,  127,
 128,  129,  130,  131,  132,  133,  134,  135,  136,  137,
 138,  139,  140,  141,  142,  143,  144,  145,  146,  147,
 148,  149,  150,  151,  152,  153,  154,  155,  156,  157 ]

slots = [
    [0,  0,  6,  0,  1],
    [1,  4,  1,  7,  0],
    [2,  3,  3,  3,  4],
    [3,  5,  3,  6, 10],
    [4,  1,  9,  7,  7],
    [5,  2,  7,  6,  9],
    [6,  2,  0,  5,  4],
    [7,  1, 10,  4,  6],
]

def makeorbits(g):

    f = [False for x in g]

    li = []
    for i in range(190):
        if f[i]:
            continue;
        j = i
        cycle = []
        while not f[j] and j <= 189:
            cycle.append(j)
            f[j] = True
            j = g[j]
        li.append(cycle)
        cycle = []
    return li

def powers(g):    # computes the powers of the group generator

    pl = []
    h = [ i for i in range(190)]
    pl.append(h)
    for iter in range(12):
        h = [h[g[i]] for i in range(190)]
        pl.append(h)
    return pl

def addtrees(fd):   # adds the two trees to the graph

    for i in [1,2,3]:
        print(f'{0:3d} {i:3d}', file=fd)
        print(f'{0+H:3d} {i+H:3d}', file=fd)
    v = 4
    for u in range(1,94):
        print(f'{u:3d} {v:3d}', file=fd)
        print(f'{u+H:3d} {v+H:3d}', file=fd)
        v += 1
        print(f'{u:3d} {v:3d}', file=fd)
        print(f'{u+H:3d} {v+H:3d}', file=fd)
        v += 1

def addslots(fd, slot):   # adds the edges that join the two trees
    global parts, powerlist

    apart = slot[0]
    u = parts[apart][0]

# slot one for this A-orbit

    bpart = slot[1]
    v = parts[bpart][0]
    voltage = powerlist[slot[2]]
    for i in range(12):
        x = powerlist[i][u]
        y = powerlist[i][voltage[v]]
        print(f'{x:3d} {y+H:3d}', file=fd)

# slot two for this A-orbit

    bpart = slot[3]
    v = parts[bpart][0]
    voltage = powerlist[slot[4]]
    for i in range(12):
        x = powerlist[i][u]
        y = powerlist[i][voltage[v]]
        print(f'{x:3d} {y+H:3d}', file=fd)

if __name__ == "__main__":

    orbits = makeorbits(g0)
    for i in range(len(orbits)):
        print(orbits[i])

    parts = orbits[-8:]
    powerlist = powers(g0)

    fd = open('g380.elist','w')
    addtrees(fd)
    for i in range(8):
        addslots(fd,slots[i])

\end{verbatim}

\section{Program 3}

The third program does the excision described above.  It reads the $(3,14)$-graph of order 380 and
outputs the $(3,14)$-graph of order 372.

\begin{verbatim}

# A Python program that does the excision step -- the final step in the
# construction of the (3,14)-graph of order 372.

def readgraph(filename):

    eli = []
    with open(filename,"r") as fd:
        for line in fd:
            a = list(map(int,line.strip().split()))
            eli.append((a[0], a[1]))
    return eli

if __name__ == "__main__":

    exset = set([0,1,2,3,190,191,192,193])
    addset = set([(4,197),(5,196),(6,199),(7,198),(8,195),(9,194)])

    elist = readgraph("g380.elist")
    newlist = [e for e in elist if e[0] not in exset and e[1] not in exset]
    for e in addset:
        newlist.append(e)

    vlist = list(range(4,190)) + list(range(194,380))
    name = {}
    for i,v in enumerate(vlist):
        name[v] = i
    with open('g372.elist','w') as fd:
        for (u,v) in newlist:
            print(f'{name[u]:3d}  {name[v]:3d}', file=fd)

\end{verbatim}


\begin{thebibliography}{9}

\bibitem{exoo14}
G. Exoo,
A small trivalent graph of girth 14,
{\em Electron. J. Combin.} 9 (2002), \#N3.

\bibitem{exjaj}
G. Exoo and R. Jajcay,
Dynamic cage survey,
{\em Electron. J. Combin.} (2013), Dynamic Survey DS16.

\bibitem{ekjs}
G. Exoo, T. Kolokolnikov, J. Janssen and T. Salamon,
Attainable bounds for algebraic connectivity and maximally connected regular graphs,
{\em J. Graph Theory} 107 (2024), 522--549.

\bibitem{nauty}
B. D. McKay and A. Piperno,
Practical graph isomorphism, II,
{\em J. Symbolic Comput.} 60 (2014), 94--112.

\bibitem{hog}
K. Coolsaet, S. D'hondt and J. Goedgebeur,
House of Graphs 2.0: A database of interesting graphs and more,
{\em Discrete Appl. Math.} 325 (2023), 97--107.
\texttt{https://houseofgraphs.org}.

\bibitem{ekhunt}
G. Exoo and T. Kolokolnikov,
The hunt for diameter-maximal cubic graphs,
in preparation.

\end{thebibliography}
\end{document}